\documentclass[12pt]{article}
\usepackage[cp1251]{inputenc}
\usepackage[english]{babel}
\usepackage{amsfonts}

\begin{document}

\def\C{{\bf C}}
\def\R{{\bf R}}
\def\N{{\bf N}}
\def\‘P{{\bf CP}}
\def\i{{\sqrt {-1}}}
\def\O{{\cal O}}
\def\A{{\cal A}}
\def\M{{\widehat{M}}}
\def\DW{{\widehat{\cal D}}}
\def\D{{\cal D}}
\def\Z{{\bf Z}}
\def\L{{\cal L}}
\def\A{{\cal A}}
\def\B{{\cal B}}
\def\P{{\cal P}}
\def\Pic{\rm{Pic}}
\def\Mat{\rm{Mat}}
\def\ord{\rm{ord}}
\def\Im{\rm{Im}}

\title{On Nonlinear Equations Integrable in Theta Functions
of Nonprincipally Polarized Abelian Varieties}
\author{A.E. Mironov
\thanks{Institute of Mathematics, 630090 Russia, Novosibirsk; e-mail:
mironov@math.nsc.ru}}
\date{}
\maketitle

\begin{abstract}
The formula of expanding the Abel variety theta function
restricted to Abel subvariety into theta functions of
this subvariety is obtained. With the help of this formula
the solution of differential equations with Jacobi
theta functions, restricted on a nonprincipally
polarized Abel subvariety and their translations
are rewritten
in terms of the theta functions
of these subvarieties.
This is exemplified by the CKP equations,
the Bogoyavlenskij system, and the Toda $g_2^{(1)}$-chain.
\end{abstract}

\section{Introduction}
In this article we show how to find solutions
of nonlinear equations in terms of theta functions
of nonprincipally polarized Abelian varieties.

Solutions of many known integrable equations are expressed
in terms of theta functions of Riemann surfaces
connected with the equations
(the Korteweg--de Vries equation, the Kadomtsev--Petviashvili equations,
geodesics on an~ellipsoid, Kovalevskaya's top,
etc.; see the survey~[1]). In some problems where
the Riemann surface admits a~holomorphic involution,
the arguments of theta functions belong to the Prym variety,
a~subvariety of the Jacobi variety of the Riemann surface,
and solutions are therefore expressed in terms of
the theta functions of the Prym variety. The Prym theta formulas
are convenient because the theta functions depend on fewer
variables and their qualitative analysis is simpler.
The Veselov--Novikov equations
and the Landau--Lifshits equations are among such equations.
In the derivation of Prym theta functional formulas
for their solutions,  the fact is essentially used that
in these cases the Prym variety is principally
polarized and therefore has, roughly speaking,
a~unique theta function.

At the same time, there are well-known equations for which the whole
dynamics reduces to nonprincipally polarized Prym subvarieties
with the polarization type
$(1,\dots,1,2,\dots,2)$.
For instance, so are the CKP equations
(the Kadomtsev--Petviashvili hierarchy of type C)~[2],
the equations of motion of a~rigid body around
a~fixed point in the Newton field with an~arbitrary (homogeneous)
quadratic potential which were integrated
by O.~I. Bogoyavlenskij~[3],
some generalized Toda chains~[4],
and the geodesic flows on quadrics and
$SO(4)$~[5,\,6].
Also, examples are known  of integrable systems that can be linearized
on Abelian varieties with the polarization types
$(1,3)$
and
$(1,6)$ (the Toda $g_2^{(1)}$-chain
and the geodesic flow on~
$SO(4)$
with a~special metric; see the survey~[7]).

In Theorem~1, we show that
for the space of theta functions
on a~nonprincipally polarized Abelian variety we can choose a~basis
that consists of lifts of theta functions with characteristics
from an~isogenous principally polarized Abelian variety.
By Theorem~2, formulas in
the Jacobi theta functions restricted to Abelian subvarieties
and their translations can be rewritten in terms of the theta functions
of the subvarieties themselves.
This is exemplified by the CKP equations,
the Bogoyavlenskij system, and the Toda $g_2^{(1)}$-chain.

The author is grateful to I.~A. Taimanov who posed
the problem for useful discussions and apt remarks.

\section{Theta Functions of Nonprincipally Polarized
Abelian Varieties}

Let
$M={\C}^g/\Lambda$
be an~Abelian variety with a~Hodge form
$\widetilde{\omega}$.
Here
$\Lambda$
is a~lattice in~
${\C}^g$.
There are real coordinates
$x_1,\dots,x_{2g}$ in ${\C}^g$
and a~form
$\omega$
cohomologous to
$\widetilde{\omega}$
of the shape
$$
\sum\limits_{j=1}^g\delta_j\,dx_j\land dx_{g+j},
$$
where
$\delta_j$
are natural numbers and
$\delta_j$
divides
$\delta_{j+1}$.
In the complex basis corresponding to the coordinates
$\delta_1x_1,\dots,\delta_gx_g$,
the lattice $\Lambda$ is given by
$\Delta{\Z}^g+\Omega{\Z}^g$,
where
$\Delta$
is the diagonal integer matrix with diagonal
$(\delta_1,\dots,\delta_g)$ and
$\Omega$
is a~symmetric
$(g\times g)$-matrix with
$\Im\Omega>0$.
The tuple
$(\delta_1,\dots,\delta_g)$
is referred to as  the {\it polarization type\/} of~$M$
and is an~invariant of the  cohomology class of the form~$\omega$.
Denote by
${\cal L}_M=\{\theta(z|\Lambda)\}$
the space of theta functions of~$M$
with the following periodicity properties:
$$
\theta(z+\lambda|\Lambda)=\theta(z|\Lambda),\quad
\lambda\in\Delta{\Z}^g,
\eqno{(1)}
$$
$$
\theta(z+\Omega e_j|\Lambda)
=\exp(-\pi i\Omega_{jj}-2\pi iz_j)\theta(z|\Lambda),
  \eqno{(2)}
$$
where
$e^{\top}_j=(0,\dots,1,\dots,0)$
(unity occupies  the $j$th position).
The dimension of
${\cal L}_M$
equals
$\delta_1\cdot\dots\cdot\delta_g$
(see, for instance,~[8,\,9]).

An~Abelian variety is {\it principally polarized\/} if
$\Delta$
is the identity matrix.
The theta function with characteristic $[a, b]$
of a~principally polarized Abelian  variety
$\widetilde{M}={\C}^g/\{{\Z}^g+\Omega{\Z}^g\}$
is defined by the absolutely convergent series
$$
\theta[a,b](z|\Omega)=\sum\limits_{n\in {\Z}^g}\exp
(\pi i \langle (n+a),\Omega(n+a) \rangle
+2\pi i
\langle (n+a),(z+b)\rangle ),
$$
where
$a, b\in {\C}^g$.
The function
$\theta[a, b](z|\Omega)$
possesses the following periodicity properties:
$$
\theta[a,b](z+m|\Omega)=\exp(2\pi i
\langle a,m\rangle)
\theta[a,b](z|\Omega),
$$
$$ \theta[a,b](z+\Omega m|
\Omega)=\exp(-2\pi i
\langle b,m \rangle
-\pi i \langle m,\Omega m\rangle -2\pi i
\langle m,z \rangle)\theta[a,b](z|\Omega),
$$
where
$m\in {\Z}^g$.
If the characteristic
$[a,b]$
is rational then the theta function with this characteristic
determines a~section of the line bundle over~$\widetilde{M}$.

{\bf Theorem~1.}
{\sl
The theta functions
$\theta[\Delta^{-1}\varepsilon, 0](z|\Omega)$,
with
$\varepsilon\in{\Z}^g/\Delta{\Z}^g$,
constitute a~basis for the space~
${\cal L}_M$.
}

{\sc Proof.}
It follows from~(1) that the theta function
$\theta(z|\Lambda)\in{\cal L}_M$
has the Fourier series expansion
$$
\theta(z|\Lambda)=\sum\limits_{m\in{\Z}^g}a_m\exp(2\pi i\langle
m,\Delta^{-1}z\rangle).
$$
Property (2) implies the following recurrent
relation for the coefficients of the series:
$$
 a_{m+\Delta e_j}= \exp(\pi i\Omega_{jj}+2\pi i\langle
m,\Delta^{-1}\Omega e_j\rangle)a_m.
\eqno{(3)}
$$
From (3) we infer that the theta function is determined by the coefficients
$\{a_m\}$, where $0\leq m_s <\delta_s$, $1\leq s\leq g$.
Denote by
$\theta_{\varepsilon}(z|\Lambda)$
the function that is given by the absolutely convergent series
$$
\sum\limits_{m\in
{\Z}^g}a_{\varepsilon+\Delta m} \exp(2\pi i \langle
\varepsilon+\Delta m,\Delta^{-1}z\rangle),
\eqno{(4)}
$$
where
$
a_{\varepsilon}=\exp(\pi i \langle \Delta^{-1}\varepsilon,
\Omega\Delta^{-1}\varepsilon\rangle)
$
and
$a_{\varepsilon+\Delta m}$
are found from the recurrent formulas~(3).  The functions
$\theta_{\varepsilon}$
constitute a~basis for ~${\cal L}_M$.
The recurrent system (3) is solvable explicitly:
$$
a_{\varepsilon+\Delta m}=
\exp(\pi i\langle m,\Omega m\rangle+
2\pi i\langle \Delta^{-1}\varepsilon,\Omega m\rangle+
\pi i\langle \Delta^{-1}\varepsilon,\Omega \Delta^{-1}\varepsilon
\rangle).
$$
Then
$$
\theta_{\varepsilon}(z|\Lambda)=
\sum\limits_{m\in {\Z}^{g}}\exp
(\pi i \langle  (m+\Delta^{-1}\varepsilon),
\Omega (m+\Delta^{-1}\varepsilon) \rangle
+2\pi i
\langle (m+\Delta^{-1}\varepsilon),z\rangle)
$$
$$
=\theta[\Delta^{-1}\varepsilon, 0](z|\Omega).
$$
The theorem is proven.

The theta functions
$\theta_{\varepsilon}=
\theta[\Delta^{-1}\varepsilon, 0](z|\Omega)$
are lifts of the theta functions
with characteristics from the principally polarized Abelian variety
$\widetilde{M}$
under the isogeny
$\xi:M\rightarrow\widetilde{M}$,\ $\xi(z)=z$,
of degree
$\delta_1\cdot\dots\cdot\delta_g$.
Observe that our construction depends on the choice of the principally
polarized Abelian variety~$\widetilde{M}$
isogenous to~$M$.
If we choose another isogeny then the lifts
of the corresponding theta functions with characteristics
as well determine a~basis for ~${\cal L}_M$.

We now find expansion for the restriction of a~theta function
of an~Abelian variety to an~Abelian subvariety
in the theta functions of the subvariety.
The curvature form of the bundle associated with the positive
divisor is the Hodge form; consequently, every
positive divisor on an~Abelian variety determines the polarization on it.

Denote by
$\widetilde{M}={\C}^g/\{\widetilde{\Delta}{\Z}^g+
\widetilde{\Omega}{\Z}^g\}$
an~Abelian variety and by
$M\subset \widetilde{M}$,
an~Abelian subvariety.
Suppose that the intersection of the theta divisor of~
$\widetilde{M}$
(the set of zeros of the theta function
$\theta(\cdot|\widetilde{\Omega})$)
with~$M$ determines the polarization type
$(\delta_1,\dots,\delta_n)$
on~$M$.
Then the restriction of~
$\theta(\cdot|\widetilde{\Omega})$
to~$M$
is a~theta function of~$M$.
Hence, there is an~isomorphism
$\varphi:{\C}^n/\{\Delta{\Z}^n+\Omega{\Z}^n\}\rightarrow M$,
where
$\Delta$
is the diagonal matrix with the diagonal
$(\delta_1,\dots,\delta_n)$
and
$\Omega$
is some symmetric matrix with
$\Im\Omega>0$,
such that
$\theta(\varphi(z)|\widetilde{\Omega})$
is a~theta function of
${\C}^n/\{\Delta{\Z}^n+\Omega{\Z}^n\}$.
Suppose that
$\varphi(z)=\Phi z$,
where
$\Phi$
is some
$(n\times g)$-matrix and
$z^{\top}=(z_1,\dots,z_n)$.
Since
$\theta(\varphi(z)|\widetilde{\Omega})$
is a~theta function of
${\C}^n/\{\Delta{\Z}^n+\Omega{\Z}^n\}$,
we have the inclusion
$\Phi\Delta\subset\widetilde{\Delta}{\Z}^g$
and the equality
$\Phi\Omega=\widetilde{\Omega}P$,
where
$P$
is some integer $(n\times g)$-matrix.
Since
$$
\theta(\varphi(z+\Omega e_j)|\widetilde{\Omega})=
\theta(\varphi(z)+\widetilde{\Omega}Pe_j|\Omega)
$$
$$
=\exp(-\pi i\langle e_j,P^{\top}\widetilde{\Omega}Pe_j\rangle-
2\pi i\langle e_j,P^{\top}\Phi z\rangle)
\theta(\varphi(z)|\widetilde{\Omega}),
$$
$P^{\top}\Phi$~
is the identity $(n\times n)$-matrix and consequently
$\Omega=P^{\top}\widetilde{\Omega}P$.

{\bf Theorem~2.}
{\sl
The following formula is valid:
$$
\theta(\varphi(z)-\gamma|\widetilde{\Omega})=
\sum\limits_{\varepsilon\in{\Z}^n/\Delta{\Z}^n}
c_{\varepsilon}
\theta[\Delta^{-1}\varepsilon,0](z-P^{\top}\gamma|\Omega),
$$
where
$\gamma^{\top}=(\gamma_1,\dots,\gamma_g)$ and
$$
c_\varepsilon
=\sum\limits_{m\in{\Z}^g|\Phi^{\top}m=\Delta^{-1}\varepsilon}
\exp(\pi i\langle m,\widetilde{\Omega} m\rangle-2\pi i
\langle m,\gamma\rangle
$$
$$
+2\pi i\langle \varepsilon, \Delta^{-1}P^{\top}\gamma\rangle-
\pi i\langle \Delta^{-1}\varepsilon,\Omega
\Delta^{-1}\varepsilon \rangle).
$$
}

{\sc Proof.}
In view of
$$
\theta(\varphi(z+\lambda)-\gamma|\widetilde{\Omega})=
\theta(\varphi(z)-\gamma|\widetilde{\Omega}),
\quad \lambda\in\Delta{\Z}^n,
$$
$$
\theta(\varphi(z+\Omega e_j)-\gamma|\widetilde{\Omega})=
\theta(\varphi(z)+\widetilde{\Omega} Pe_j-\gamma|\widetilde{\Omega})
$$
$$
=\exp(-\pi i\langle Pe_j,\widetilde{\Omega} Pe_j\rangle-
2\pi i\langle Pe_j,\Phi z-\gamma\rangle)
\theta(\varphi(z)-\gamma|\widetilde{\Omega})
$$
$$
=\exp(-\pi i\langle e_j,\Omega e_j\rangle-
2\pi i\langle e_j, z-P^{\top}\gamma\rangle)
\theta(\varphi(z)-\gamma|\widetilde{\Omega}),
$$
the function
$\theta(\varphi(z)-\gamma|\widetilde{\Omega})$
is a~theta function of the Abelian variety~$M$
with argument
$z-P^{\top}\gamma$
and hence it expands in the basis theta functions~(4).
The  Fourier series expansion of the function
$\theta(\varphi(z)-\gamma|\widetilde{\Omega})$
looks like
$$
\sum\limits_{m\in{\Z}^g}
\exp(\pi i\langle m,\widetilde{\Omega} m\rangle+
2\pi i
\langle m,\Phi z-\gamma\rangle).
$$
By (4), the coefficient of
$
\exp(2\pi i\langle\Delta^{-1}\varepsilon,
z-P^{\top}\gamma\rangle)
$
in this expansion equals
$
c_{\varepsilon}\exp(\pi i \langle \Delta^{-1}\varepsilon,
\Omega\Delta^{-1}\varepsilon\rangle)
$.
Whence we obtain the formula for
$c_{\varepsilon}$.
The theorem is proven.

\section{The Theorem on Expansion of the Prym Theta
Function}

Here we confine exposition to the case in which
the nonprincipally polarized Abelian variety~$M$
is a~Prym variety.

Suppose that
$\Gamma$
is a~Riemann surface with an~involution
$\sigma:\Gamma\rightarrow \Gamma$
having
$2(n+1)$
fixed points
$Q_0,\dots,Q_{2n+1}$, $n>0$,
and
$\pi:\Gamma\rightarrow\Gamma_0=\Gamma/\sigma$
is a~projection.
The genus of~
$\Gamma$
equals
$2g+n$,
where
$g$
is the genus of~$\Gamma_0$.
The surface
$\Gamma$
possesses a~canonical basis of cycles
$$
a_1,\dots,a_{g+n},
\tilde{a}_1,\dots,\tilde{a}_g,
b_1,\dots,b_{g+n},
\tilde{b}_1,\dots,\tilde{b}_g
$$
such that
$$
\sigma(a_{\alpha})+\tilde{a}_{\alpha}=
\sigma(b_{\alpha})+\tilde{b}_{\alpha}=0,
\quad
1\leq\alpha\leq g,
$$
$$
\sigma(a_j)+a_j=
\sigma(b_j)+b_j=0,
\quad
g+1\leq j\leq g+n.
$$
To the basis of cycles, there corresponds the canonical
basis of Abelian differentials
$$
u_1,\dots,u_{g+n},\tilde{u}_1,\dots,\tilde{u}_g
$$
satisfying the relations
$$
\sigma^*u_{\alpha}+\tilde{u}_{\alpha}=0,\quad
\sigma^*u_j+u_j=0,
\quad 1\leq\alpha\leq g,\quad
g+1\leq j\leq g+n.
$$
The cycles
$a_1,\dots,a_g,b_1,\dots,b_g$
project to the canonical basis of cycles on~$\Gamma_0$.
The surface
$\Gamma_0$
possesses the canonical basis of Abelian differentials
$\widetilde{\omega}_1,\dots,\widetilde{\omega}_g$
such that
$\pi^*(\widetilde{\omega}_k)=u_k-\tilde{u}_k$, $1\leq k\leq g$.
Denote by~
$T$
the matrix of periods of~
$\widetilde{\omega}_k$,
$T_{jk}=\int\nolimits_{\pi(b_j)}\widetilde{\omega}_k$,
$1\leq j$, $k\leq g$.
Recall that
$\int\nolimits_{\pi(a_j)}\widetilde{\omega}_k=\delta_{jk}$
by the choice of the basis cycles and the basis Abelian
differentials.
Denote by
$J={\C}^{2g+n}/\{{\Z}^{2g+n}+\Omega{\Z}^{2g+n}\}$
the Jacobi variety of ~$\Gamma$,
where
$\Omega$
is the matrix of periods of the basis
differentials of ~$\Gamma$,
and denote by
$A:\Gamma\rightarrow  J$
the Abelian mapping with basepoint~$Q_0$.
The differentials
$$
\omega_{\alpha}=u_{\alpha}+\tilde{u}_{\alpha},
\quad
\omega_j=u_j,\quad
1\leq\alpha\leq g,\quad g+1\leq j\leq g+n,
$$
constitute a~basis for the Abelian Prym differentials
$\sigma^*\omega_k=-\omega_k$,\ $1\leq k\leq g+n$.
The involution
$\sigma$
indices the involution
$\sigma_*:J\rightarrow J$,
$$
\sigma_*(z_1,\dots,z_g,z_{g+1},\dots,z_{g+n},
\tilde{z}_1,\dots,\tilde{z}_g)=
-(\tilde{z}_1,\dots,\tilde{z}_g,z_{g+1},\dots,z_{g+n},
z_1,\dots,z_g).
$$
The Prym variety is the Abelian subvariety
$$
Pr=\{z\in J\ |\ \sigma_*(z)=-z\}\subset J.
$$
Denote by
$\varphi:{\C}^{g+n}/\Lambda \rightarrow Pr$
the isomorphism
$$
\varphi(z_1,\dots,z_{g+n})=
({1}/{2}z_1,\dots,{1}/{2}z_g,z_{g+1},\dots,z_{g+n},
{1}/{2}z_1,\dots,{1}/{2}z_{g+n}),
$$
where
$
\Lambda=\Delta{\Z}^{g+n}+\Pi{\Z}^{g+n}$,
$\Delta$
is the diagonal matrix with diagonal
$(2,\dots,2$ , $1, \dots,1)$
($n$ units), and
$\Pi$
is the symmetric matrix with
$\Im \Pi>0$
whose entries are
$$
\Pi_{\alpha k}=2\int\limits_{b_\alpha}\omega_k,
\quad 1\leq\alpha\leq g,
\quad 1\leq k\leq g+n,
$$
$$
\Pi_{jk}=\int\limits_{b_j}\omega_k,
\quad g+1\leq j\leq g+n,
\quad
1\leq k\leq g+n.
$$
The variety
$Pr$
is nonprincipally polarized and the dimension of the space
of the Prym theta functions equals~$2^g$.

From Theorems~1 and 2 we derive the following

{\bf Theorem~3.}
{\sl
The theta functions
$\theta[\Delta^{-1}\varepsilon, 0](z|\Pi)$,
where
$\varepsilon=(\varepsilon_1,\dots,\varepsilon_g$, $0,\dots,0)$
$(n$ zeros$)$,
$\varepsilon_j\in\{0,1\}$,
constitute a~basis for the space of the Prym theta functions.
The following expansion in the Prym theta functions is valid:
$$
\theta(\varphi(z)-\gamma|\Omega)=
\sum\limits_{\varepsilon}c_{\varepsilon}
\theta[\Delta^{-1}\varepsilon, 0]
(z-\tilde{\gamma}|\Pi),
\eqno{(5)}
$$
where
$ \gamma=(\gamma_1,\dots,\gamma_{g+n},
\tilde{\gamma}_1,\dots,\tilde{\gamma}_g)$,
$\tilde{\gamma}=(\gamma_1+\tilde{\gamma}_1,\dots,
\gamma_g+\tilde{\gamma}_g,\gamma_{g+1}\dots,\gamma_{g+n})$,
$$
c_{\varepsilon}=
\sum\limits_{m\in{\Z}^g}
\exp(\pi i\langle m_{\varepsilon},\Omega m_{\varepsilon}\rangle
+\pi i\langle\varepsilon, \tilde{\gamma}\rangle-
2\pi i\langle m_{\varepsilon},\gamma\rangle
-\pi i
\langle \Delta^{-1}\varepsilon, \Pi\Delta^{-1}\varepsilon\rangle),
\eqno{(6)}
$$
and
$m_{\varepsilon}$
stands for the vector
$(m_1,\dots,m_g,0,\dots,0,\varepsilon_1-m_1,\dots,\varepsilon_g-m_g)$
$(n$ zeros$)$.
}

In~[10, Proposition 5.5], some formula was obtained that connects
theta functions with the characteristics
of the principally polarized Abelian varieties
$J$,
${\C}^g/\{{\Z}^g+T{\Z}^g\}$,
and
${\C}^{g+n}/\{{\Z}^{g+n}+\Pi{\Z}^{g+n}\}$.
This formula gives another way to derivation
of ~(5) and implies that the constants~
$c_{\varepsilon}$
(in Theorem~3) equal
$\theta[\tilde{\varepsilon}, 0](\hat{\delta}| 2T)$,
where
$\tilde{\gamma}=(\tilde{\gamma}_1-\gamma_1,\dots,
\tilde{\gamma}_g-\gamma_g)$
and
$\tilde{\varepsilon}=\frac{1}{2}(\varepsilon_1,\dots,\varepsilon_g)$.

\section{Applications}

{\bf 4.1.\ The CKP hierarchy.}
This hierarchy is determined by an~infinite system of Lax equations
$$
[\partial_n-B_n,\partial_m-B_m]=0,\quad m,n=1,3,5,\dots,
$$
in the coefficients of the operators
$$
B_n=\partial^n+\sum\limits_{i=0}^{n-2}u_{ni}\partial^i
$$
that depend on infinitely many variables
$x=t_1, t_3, t_5,\dots$;
moreover, the equality
$B_n^*=-B_n$
must hold. Here  $B_n^*$
is the formal adjoint of~$B_n$,
$\partial={\partial}/{\partial x}$,
and
$\partial_n={\partial}/{\partial t_n}$.
There is a~pseudodifferential operator
$$
L=\partial+\frac{V}{2}\partial^{-1}+\sum\limits_{k=2}^{\infty}V_k
\partial^{-k}
$$
such that
$B_k=(L^k)^+$,
where
$(L^k)^+$
is the differential part of~$L^k$ and
the functions
$V_k$
are expressed in terms of
$V$
and its derivatives.
The first two operators of the hierarchy are
$$
B_3=\partial^3+\frac{3}{2}V\partial+\frac{3}{4}\partial V,
\quad
B_5=\partial^5+\frac{5}{2}V\partial^3+
\frac{15}{4}\partial V\partial^2+
W\partial+\frac{\partial W}{2}+\frac{5}{8}\partial^3V,
$$
where
$$
W=\frac{1}{3\partial^2V}\left(\frac{\partial^6 V}{3}+
5V\partial^4V+\frac{45}{4}\partial V\partial^3V+
15(\partial^2 V)^2+\frac{15}{2}V(\partial V)^2 \right.
$$
$$
\left.+\frac{15}{2}V^2\partial^2 V-
\frac{5}{3}\partial^2_3V
-\frac{5}{3}\partial_3\partial^3V-
5\partial_3V\partial V-\frac{5}{2}V\partial_3\partial V+
3\partial_5\partial V\right).
$$
The first equation of the hierarchy is
$$
\partial_3V=\frac{6}{5}\partial W-\frac{7}{2}\partial^3 V-
3V\partial V.
$$
Express the function
$V$
(a~solution to the CKP hierarchy) in terms of the theta functions
of the Prym variety.
Let
$\Gamma$
be a~Riemann surface with an~involution
$\sigma:\Gamma\rightarrow \Gamma$ having
$2(n+1)$
fixed points
$Q_0,\dots,Q_{2n+1}$, $n>0$.
Fix a~nonspecial divisor
$D=P_1+\dots +P_{2g+n}$
on~
$\Gamma$
and a~polynomial~$R$.
Take a~local parameter
$k^{-1}$
in a~neighborhood of ~$Q_1$ such that
$k\sigma=-k$.
The one-point Baker--Akhiezer function with spectral data
$\{\Gamma,Q_1,D,k^{-1},R(k)\}$
is a~function
$\psi(P)$, $P\in\Gamma$,
defined (to within proportionality)
by the following properties~[11]:

(1)\ $\psi(P)$
is meromorphic on
$\Gamma\backslash Q_1$
and the set of its poles coincides with~$D$;

(2)\ the function
$\psi(P)\exp(-R(k))$
is analytic in a~neighborhood of~$Q_1$.

{\bf Theorem A~\rm[2]}
{\sl
If a~nonspecial positive divisor~$D$ satisfies the relation
$$
D+\sigma D \sim C_{\Gamma}+2Q_1,
$$
where
$C_{\Gamma}$
is the canonical class on~$\Gamma$,
then the one-point Baker--Akhiezer function
constructed from the spectral data
$\{\Gamma,Q_1,D,k^{-1},xk+t_3k^3+t_5k^5+\dots\}$
is an~eigenfunction of the operators of the {\rm CKP} hierarchy;
i.e.,
$$
B_m\psi=\partial_m\psi.
$$
}

Recall that
$C_{\Gamma}$
is the linear equivalence class of the divisor
of zeros and poles of some meromorphic 1-form. The condition
of the theorem means that there is a~meromorphic 1-form~$\omega_1$
with zeros in~$D+\sigma D$
and with a~second-order pole at~$Q_1$.

Let
$\omega_s$
be a~Prym differential of the second kind with a~unique pole
of order
$s+1$
($s$
is an~odd number) at ~$Q_1$
and
with the zero~$a$- and
$\tilde{a}$-periods.
Denote by
$U_s$
the vector
$(U_{s1},\dots,U_{sg},U_{s\,g+1},
\dots,U_{s\,g+n},  U_{s1},\dots,U_{sg})$,
where
$U_{s_j}=\int\nolimits_{b_j}\omega_s$.
The Baker--Akhiezer function
$\psi(x,t_3,t_5,\dots;P)$
looks like~[11]
$$
\exp\biggl(2\pi ix\int\limits_{Q_0}^P\omega_1+2\pi
it_3\int\limits_{Q_0}^P\omega_3+ 2\pi
it_5\int\limits_{Q_0}^P\omega_5+\dots\biggr)
$$
$$
\times\frac{\theta(A(P)-A(D)-K_{\Gamma}+xU_1+t_3U_3+t_5U_5+\dots|\Omega)}
{\theta(A(P)-A(D)-K_{\Gamma}|\Omega)},
$$
where
$K_{\Gamma}=-\frac{1}{2}A(C_{\Gamma})$
is the vector of the Riemann constants and
$\theta(.|\Omega)$
is the theta function of the Jacobi variety of ~$\Gamma$.
From the  power series expansion of ~
$\psi$
in a~neighborhood of
~$Q_1$ we obtain the following

{\bf Corollary 1 \rm (to Theorem~A)}
{\sl
The finite-gap solutions to~{\rm CKP} look like
$$
V(x,t_3,t_5,\dots)=
2\partial^2\log\theta(xU_1+t_3U_3+t_5U_5+\dots-\gamma|\Omega),
$$
where
$\gamma=A(D)+K_{\Gamma}$.
}

From Theorem~3 we derive the next

{\bf Theorem~4}
{\sl
The finite-gap solutions to~{\rm CKP}
are expressed in terms of the Prym theta functions by the formula
$$
V(x,t_3,t_5,\dots)=
2\partial^2\log
\sum\limits_{\varepsilon}c_{\varepsilon}
\theta[\Delta^{-1}\varepsilon, 0]
(x\tilde{u}_1+
t_3\widetilde{U}_3+
t_5\widetilde{U}_5+\dots-\tilde{\gamma}|\Pi),
$$
where
$\widetilde{U}_s=
(2U_{s1},\dots,2U_{sg},U_{s\,g+1},\dots,U_{s\,g+n})$,
$\varepsilon=
(\varepsilon_1,\dots,\varepsilon_g,0,\dots,0)$
$(n$ zeros$)$,
$\varepsilon_j\in\{0,1\}$, and
$c_{\varepsilon}$
are found by~{\rm (6)}.
}

It would be interesting to explain
a~solution to CKP in terms of secants of Abelian varieties.
This was done in~[12]
for soliton equations integrable in
theta functions of Jacobi varieties
and in~[13] for equations integrable in theta
functions of principally polarized
Prym varieties (see~also~[9]).

{\bf 4.2.\ The problem of rotation of a~rigid body.}
Rotation of a~rigid body
$S$
around a~fixed point~
$O\in S$
in a~coordinate system
$r_1$, $r_2$, $r_3$
rotating together with~$S$
in the Newton field with potential~$\varphi$
is described by the generalized Euler equations~[3]:
$$
\dot{M}=M\times\omega+
\frac{\partial U}{\partial\alpha}\times\alpha+
\frac{\partial U}{\partial\beta}\times\beta+
\frac{\partial U}{\partial\gamma}\times\gamma,
 \eqno{(7)}
$$
$$
\dot{\alpha}=\alpha\times\omega,\quad
\dot{\beta}=\beta\times\omega,\quad
\dot{\gamma}=\gamma\times\omega,
$$
$$
U(\alpha,\beta,\gamma)=
\int\limits_S\rho(r)\varphi(\langle r,\alpha\rangle,
\langle r,\beta\rangle,\langle r,\gamma\rangle)\,
dr_1dr_2dr_3,
$$
where
$\rho(r)$
is the density of the body~$S$
at~
$r=(r_1, r_2, r_3)$;
$\alpha$, $\beta$, $\gamma$
is the orthonormal basis of the fixed coordinate system;
and the vectors $M$
and
$\omega$
of kinetic momentum and angular velocity
are connected by the relations
$$
M_i=\sum\limits_{k=1}^3I_{ik}\omega_k,
$$
with
$I_{ik}$
the components of the inertia tensor of the body
$S$
in the rotating coordinate system:
$$
I_{ik}=\int\limits_S\rho(r)
\Biggl(\delta_{ik}\sum\limits_{j=1}^3r_j^2-r_ir_k\Biggr)\,
d r_1 d r_2 d r_3.
$$
On using the isomorphism between the Lie algebra~
${\R}^3$
with vector product
and the algebra of skew-symmetric $(3\times 3)$-matrices
with the commutator product,  it was demonstrated in~[3]
that (7) yields the matrix equations
$$
\biggl[L,\frac{\partial}{\partial t}+Q\biggr]=0,\quad
L=BE^2+ME+u,\quad Q=\omega-EI,
$$
where
$u$ and~$B$
are some symmetric matrices and
$E$
is an~arbitrary parameter. For simplicity,
we use the same symbols to denote
the skew-symmetric matrices and the  corresponding vectors
under this isomorphism.
Denote by
$\Gamma$
the smooth completion of the surface that is defined in
${\C}^2$
by the equation
$$
\det(BE^2+ME+u-w1)=0.
$$
The Riemann surface~$\Gamma$
is not hyperelliptic and
admits the holomorphic involution
$\sigma:\Gamma\rightarrow\Gamma$, $\sigma(w,E)=(w,-E)$.
The genus of ~$\Gamma$
equals~4, the surface~
$\Gamma_0=\Gamma/\sigma$
is elliptic, and the involution~$\sigma$
has 6 fixed points. The dimension of the Prym variety equals three, and
the components~
$\omega_i^j(t)$
of the angular velocity are equal to
$$
A_i^j\exp\bigl(t\xi_i^j\bigr)\frac{\theta\bigl(tU+z_i^j|\Omega\bigr)}
{\theta(tU+z_0|\Omega)},
$$
where the constants~$A_i^j$ and $\xi_i^j$
and the vectors
$z_i^j$, $z_0\in J(\Gamma)$, and $U\in Pr(\Gamma)$
are defined in~[3].
From Theorem~3 we derive the following

{\bf Theorem~5}
{\sl
The components of the angular velocity
are expressed in terms of the
Prym theta functions by the formula
$$
\omega_i^j(t)= A_i^j\exp\bigl(t\xi_i^j\bigr)
\frac{c_0
\theta\bigl(t\widetilde{U}+\tilde{z}_i^j|\Pi\bigr)+
c_1\theta[\Delta^{-1}\varepsilon_1,0]
\bigl(t\widetilde{U}+\tilde{z}_i^j|\Pi\bigr)}
{c'_0\theta(t\widetilde{U}+\tilde{z}|\Pi)
+c'_1\theta[\Delta^{-1}\varepsilon_1,0]
(t\widetilde{U}+\tilde{z}|\Pi)},
$$
where
$\varepsilon_1=(1,0,0)$,
$\Delta=(2,1,1)$,
and the constants
$c_*$ and~ $c'_*$
are defined by~{\rm (6)}.
}

{\bf\boldmath 4.3. The Toda $g_2^{(1)}$-chain.}
This chain is described by the Hamiltonian system with Hamiltonian
$$
H=\frac{1}{2}\sum\limits_{i=1}^3p_i^2+\exp(q_2-q_1)+\exp(q_3-q_2)+
\exp\frac{1}{3}(q_1+q_2-2q_3)
$$
which is reduced by a~change of coordinates to the system
$$
\dot{Y}=CX,
\quad \dot{X}^{\top}=(x_1y_1,x_2y_2,x_3y_3),
 \eqno{(8)}
$$
where
$X^{\top}(t)=(x_1(t),x_2(t),x_3(t))$,
$Y^{\top}(t)=(y_1(t),y_2(t),y_3(t))$,
and
$C$
is the Cartan matrix of the Kac--Moody Lie $g_2^{(1)}$-algebra
which equals
$$
\left(
\begin{array}{ccc}
 2  & -1 &  0 \\
 -1 & 2  & -3 \\
 0  & -1 &  2
\end{array} \right) .
$$
System (8) admits the Lax representation
$$
\dot{A_{\mu}}=[A_{\mu},B_{\mu}], \eqno{(9)}
$$
where
$$
A_{\mu}=
\left(
\begin{array}{ccccccc}
 b_1 & a_2 & \mu^{-1}a_3 & 0 & a_1 & 0 & 0 \\
 a_2 & b_2 & 0 & \sqrt{2}a_1 & 0 & 0 & 0 \\
 \mu a_3 & 0 & b_3 & 0 & 0 & 0 & -a_1 \\
 0 & \sqrt{2}a_1 & 0 & 0 & 0 & -\sqrt{2}a_1 & 0 \\
 a_1  & 0 & 0 & 0 & -b_3 & 0 & -\mu^{-1}a_3 \\
 0 & 0 & 0 & -\sqrt{2}a_1  & 0 & -b_2 & -a_2 \\
 0 & 0 & -a_1 & 0 & -\mu a_3 & -a_2 & -b_1
\end{array}\right),
$$
$$
B_{\mu}=
\left(
\begin{array}{ccccccc}
 0 & a_2 & -\mu^{-1}a_3 & 0 & -a_1 & 0 & 0 \\
 -a_2 & 0 & 0 & \sqrt{2}a_1 & 0 & 0 & 0 \\
 \mu a_3 & 0 & 0 & 0 & 0 & 0 & a_1 \\
 0 & -\sqrt{2}a_1 & 0 & 0 & 0 & -\sqrt{2}a_1 & 0 \\
 a_1  & 0 & 0 & 0 & 0 & 0 & \mu^{-1}a_3 \\
 0 & 0 & 0 & \sqrt{2}a_1  & 0 & 0 & -a_2 \\
 0 & 0 & -a_1 & 0 & -\mu a_3 & a_2 & 0
\end{array}\right),
$$
$$
a_1=\frac{i}{2}\sqrt{x_3},\quad
a_2=\frac{i}{2}\sqrt{x_2},\quad
a_3=\frac{i}{2}\sqrt{x_1},
$$
$$
b_1=\frac{y_1+y_3}{4},\quad
b_2=\frac{y_1-2y_2+y_3}{4},\quad
b_3=\frac{3y_1+y_3}{4}.
$$
It follows from (9) that the operators~
$A_{\mu}$
and
$\partial_t+B_{\mu}$
commute; consequently, the operator~$A_{\mu}$
takes the seven-dimensional kernel of the operator
$\partial_t+B_{\mu}$
into itself.
Hence, the eigenvalues of ~$A_{\mu}$
are independent of time. Therefore, the characteristic polynomial
$Q(\mu,\lambda)=\det(A_{\mu}-\lambda E)$
of ~$A_{\mu}$
is independent of time:
$$
Q(\mu,\lambda)=\lambda
(H_1(X,Y)({1}/{\mu}+\mu)-
\lambda^6-H_2(X,Y)\lambda^4-H_3(X,Y)\lambda^2-H_4(X,Y)).
$$
The coefficients
$H_i(X,Y)=c_i$
(functions of the components
$A_{\mu}$)
are integrals of motion. Denote by~$\Gamma$
the smooth completion of the curve given in~${\C}^2$
by the equation
$$
c_1({1}/{\mu}+\mu)=
\lambda^6+c_2\lambda^4+c_3\lambda^2+c_4. \eqno{(10)}
$$
The spectral curve~$\Gamma$ admits the two involutions
$$
\tau :(\mu,\lambda)\rightarrow ({1}/{\mu},\lambda),
\quad
\sigma :(\mu,\lambda)\rightarrow (\mu,-\lambda).
$$
The involution
$\tau$
is hyperelliptic with fixed points
$P_i(1,\lambda_i)$, $1\leq i\leq 6$,
where
$\lambda_i$
is a~root of the equation
$\lambda^6+c_2\lambda^4+c_3\lambda^2+c_4-2c_1=0$,
and
$P_j(-1,\lambda_j)$, $7\leq j\leq 12$,
where
$\lambda_j$
is a~root of the equation
$\lambda^6+c_2\lambda^4+c_3\lambda^2+c_4+2c_1=0$.
The involution~$\sigma$
has~4 fixed points:
$R_1(\mu_1,0)$, $R_2(\mu_2,0)$, $R_3(0,\infty)$,
and
$R_4(\infty,\infty)$,
where
$\mu_1$
and
$\mu_2$
are roots of the equation
$c_1({1}/{\mu}+\mu)=c_4$,
$R_3$ and $R_4$
are points at infinity of the curve given by equation~(10).
Finite-gap solutions to equation~(9) are expressed in terms of
the theta function of the Jacobi variety~$J$ of the curve~$\Gamma$
(see~[14]).
The arising winding of the torus~$J$
covers a~two-dimensional
Abelian subvariety~$M$
with the polarization type (1,3) [4,\,7]
which is contained in the Prym variety of the involution~$\sigma$.
 By Theorem~2,
solutions to~(9) can consequently be expressed in terms of
the theta function of the Abelian variety~$M$.

{\bf References.}

[1]
 Dubrovin~B.~A., Krichever~I.~M., and Novikov~S.~P.,
Integrable systems.~I,
Contemporary  Problems of Mathematics. Fundamental Trends. Dynamical Systems
[in~Russian], VINITI, Moscow, 1985, {\bf4}, pp.~179--284.
(Itogi Nauki i Tekhniki.)

[2]
 Date~E., Jimbo~M., Kashiwara~M., and  Miwa~T.,
      Transformation groups for soliton equations,
       J. ~Phys. Soc. Japan,
       1981,
      vol.    50,
      3813--3818.

[3]
 Bogoyavlenskij~O.~I.,
      Integrable Euler equations on Lie algebras
arising in the problems of mathematical physics,
       Izv. Akad. Nauk SSSR Ser. Mat.,
       1984,
      vol. 48,
      N. 5,
      883--938.

[4]
 Adler~M. and  Van Moerbeke~P.,
      Linearization of Hamiltonian systems, Jacobi varieties, and
representation theory,
       Adv. Math.,
       1980,
      vol. 38,
       318--379.

[5]
 Audin ~M.,
      Courbes algebriques et systemes
integrables:  geodesiques des quadriques,
       Expo\-si\-tio\-nes Math.,
       1994,
      vol. 12,
       193--226.

[6]
 Haine L.,
      Geodesic flow on $SO(4)$ and Abelian surfaces,
       Math. Ann.,
      1983,
      vol.    263,
       435--472,

[7]
 Van Moerbeke~P.,
Introduction to algebraic integrable systems and their Painleve
analysis,
 Proceedings of the Symposium in Pure Mathematics,
1989, Vol.~49, part~1, 107--131.

[8]
 Griffiths~F. and Harris~J.,
       Principles of Algebraic Geometry,
       1982,
        Mir,
        Moscow,
    Russian translation.

[9]
 Taimanov~I.~A.,
      Secants of abelian varieties, theta-functions,
      and soliton equations,
      Uspekhi Mat. Nauk,
      1997,
      vol.    52,
      N.  1,
          149--224.

[10]
 Fay~J.~D.,
Theta functions on Riemann surfaces,
Lecture Notes in Math.; {\bf 352},
Springer-Verlag, Berlin, Heidelberg, and  New York (1973).

[11]
 Krichever~I.~M.,
       The methods of algebraical geometry
      in the theory of nonlinear equations,
       Uspekhi Mat. Nauk,
       1977,
      vol. 32,
      N. 6,
       183--208.

[12]
 Mumford~D.,
       Tata Lectures on Theta,
       1988,
        Mir,
        Moscow,
 Russian translation.

[13]
 Taimanov~I.~A.,
       Prym theta-functions and hierarchies of nonlinear equations
       Mat. Zametki,
       1991,
      vol.    50,
      N. 1,
        98--107.

[14]
 Krichever~I.~M.,
Nonlinear equations and elliptic curves,
 Contemporary  Problems of Mathematics
[in~Russian], VINITI, Moscow, 1983, {\bf23}, 79--136.

\enddocument